\documentclass[11pt]{article}

\usepackage{theorem,amstex,amscd}

\title{Integrability of the twistor space for a hypercomplex manifold.}

\author{D. Kaledin}

\date{Independent University of Moscow}

\makeatletter
\renewcommand{\subsection}{\@startsection{subsection}{2}{0pt}{-3.25ex
plus -1ex minus -.2ex}{-0pt}{\normalfont\normalsize\bfseries}}
\@addtoreset{equation}{section}
\makeatother

\newcommand{\C}{{\Bbb C}}
\newcommand{\R}{{\Bbb R}}
\newcommand{\h}{{\Bbb H}}
\newcommand{\I}{{\cal I}}
\newcommand{\E}{{\cal E}}
\newcommand{\V}{{\cal V}}
\newcommand{\A}{{\cal A}}
\newcommand{\calo}{{\cal O}}
\newcommand{\cp}{{{\Bbb C}P^1}}

\newcommand{\wt}{\widetilde}
\newcommand{\wh}{\widehat}
\newcommand{\6}{\partial}

\newcommand{\id}{\operatorname{id}}
\newcommand{\Maps}{\operatorname{Maps}}
\newcommand{\SB}{{\operatorname{SB}}}

\newtheorem{theorem}{Theorem}
\newtheorem{lemma}[theorem]{Lemma}

\theorembodyfont{\upshape\rm}

\newtheorem{defn}{Definition}[section]

\newcommand{\proof}{\noindent{\em Proof. }}
\def\endproof{\ensuremath{\square}}

\def\eqref#1{\thetag{\ref{#1}}}

\begin{document}

\maketitle

\section*{Introduction}

A {\em hyperk\"ahler manifold} is by definition a Riemannian
manifold eqipped with a smooth parallel action of the algebra of
quaternions on its tangent bundle.  Hyperk\"ahler manifolds were
introduced by Calabi in \cite{C} and have since been the subject of
much research. They have been shown to possess many remarkable
properties and a rich inner structure. (See \cite{Bes}, \cite{HKLR}
for an overview.) In particular, there is a so-called {\em twistor
space} associated to every hyperk\"ahler manifold. The twistor space
is a complex manifold equipped with some additional structures, and
many of the differential-geometric properties of a hyperk\"ahler
manifold can be described in terms of holomorphic properties of its
twistor space.

Some of the properties of hyperk\"ahler manifolds do not actually
depend on the Riemannian metric but only on the quaternion
action. In particular, for every manifold equipped with a smooth
action of quaternions (or, for brevity, a {\em quaternionic
manifold}) one can construct an almost complex manifold which becomes
the twistor space in the hyperk\"ahler case. Thus it would be very
convenient to have a notion of ``a hyperk\"ahler manifold without a
metric'', in the same sense as a complex manifold is a K\"ahler
manifold without a metric. The analogy with the K\"ahler case
suggests that this would require a certain integrability condition
on the quaternionic action, automatic in the Riemannian case. One
version of such a condition was suggested in \cite{Bes}, but the
resulting notion of an integrable quaternionic manifold is too
restrictive and excludes many interesting examples.

A more convenient notion is that of a {\em hypercomplex manifold}
(see \cite{Bo}). By definition a hypercomplex manifold is a smooth
manifold equipped with two integrable anticommuting almost complex
structures. (Note that two anticommuting almost complex structures
induce an action of the whole quaternion algebra on the tangent
bundle, and their integrability is in fact a condition on the
resulting quaternionic manifold.) In this paper we show that this
condition is in fact equivalent to the integrability of the almost
complex twistor space associated to the quaternionic manifold in
question.

Here is a brief outline of the paper. In Section~\ref{1} we give the
necessary definitions and formulate the result
(Theorem~\ref{main}). In Section~\ref{2} we describe a
linear-algebraic construction somewhat analogous to the Borel-Weyl
localization of finite-dimensional representations of a reductive
group. In Section~\ref{3} we use this version of localization to
prove Theorem~\ref{main}. The paper is essentially self-contained
and does not require any prior knowledge of the theory of
hyperk\"ahler manifolds.

\noindent
{\bf Acknowledgements.} It is a pleasure to express my gratitude to
Misha Verbitsky and Tony Pantev for many interesting and stimulating
discussions and constant encouragement. I am especially grateful to
Misha for his valuable suggestions on the present paper. 

\section{Preliminaries.}\label{1}

\subsection{}
Let $\h$ be the algebra of quaternions. 

\begin{defn} 
A {\em quaternionic manifold} is a smooth manifold $M$ equipped with
a smooth action of the algebra $\h$ on the tangent bundle
$\Theta(M)$ to $M$.
\end{defn}

Let $M$ be a quaternionic manifold. Every algebra embedding $I:\C \to \h$ 
defines by restriction an almost complex structure on $M$. Call it {\em the 
induced almost complex structure} and denote it by $M_I$. 

\subsection{}\label{maps}
The set $\Maps(\C,\h)$ of all algebra embeddings $\C \to \h$ can be
given a natural structure of a complex manifold as follows. The
algebra $\h \otimes_\R \C$ is naturally isomorphic to the $2 \times
2$-matrix algebra over $\C$. Every algebra embedding $I:\C \to \h$
defines a structure of a $2$-dimensional vector space $\h_I$ on $\h$
by means of left multiplication by $I(\C)$. It also defines a
$1$-dimensional subspace $I(\C) \subset \h_I$. The action of $\h$ on
itself by right multiplication preserves the complex structure
$\h_I$ and extends therefore to an action of the matrix algebra $\h
\otimes_\R \C$.
 
Let $\wh{I} \subset \h \otimes_\R \C$ be the annihilator of $I(\C)
\subset \h_I$. The ideal $\wh{I} \subset \h$ is a maximal right
ideal, moreover, we have $\h_I = \h \otimes \C / \wh{I}$.  It is
easy to check that every maximal right ideal of $\h \otimes_\R \C$
can be obtained in this way. This extablishes a bijection between
$\Maps(\C,\h)$ and the set of maximal right ideals in $\h \otimes_\R
\C$. It is well-known that this set coincides with the complex
projective line $\cp$.

\subsection{}
Let now $M$ be a quaternionic manifold, and let $X = M \times \cp$ be the 
product of $M$ with the smooth manifold underlying $\cp$. For every point 
$x = m \times I\in M \times \cp$ the tangent bundle $T_xX$ decomposes 
canonically as $T_xX = T_mM \oplus T_I\cp$. Define an endomorphism 
$\I:T_xX \to T_xX$ as follows: it acts as the usual complex structure map 
on $T_I\cp$, and on $T_mM$ it acts as the induced complex structure map 
$I:T_mM \to T_mM$ corresponding to $I \in \cp \cong \Maps(\C,\h)$. 
The map $\I$ obviously depends smoothly on the point $x \in X$ and 
satisfies $\I^2 = -1$. Therefore, it defines an almost complex structure on
the smooth manifold $X$. 

\begin{defn}
The almost complex manifold $X = M \times \cp$ is called {\em the twistor 
space} of the quaternionic manifold $M$. 
\end{defn}

\subsection{}
The twistor space $X$ has the following obvious properties.
\begin{enumerate}
\item The canonical projection $\pi:X \to \cp$ is compatible with
the almost complex structures. 
\item For every point $m \in M$ the embedding $\wt{m} = m \times \id:
\cp \to M \times \cp = X$ is compatible with the almost complex structures. 
\end{enumerate}
The embedding $\wt{m}:\cp \to X$ will be called {\em the
twistor line} corresponding to the point $m \in M$.

\subsection{}
The goal of this paper is to prove the following theorem. 

\begin{theorem}\label{main}
Let $M$ be a quaternionic manifold, and let $X$ be its twistor space. 
The following conditions are equivalent: 
\begin{enumerate}
\item For two algebra maps $I,J:\C \to \h$ such that $I \neq J$ and 
$\overline{I} \neq J$ the induced almost complex structures $M_I$, $M_J$ on 
$M$ are integrable. 
\item For every algebra map $I:\C \to \h$ the induced almost complex
structure $M_I$ on $M$ is integrable. 
\item The almost complex structure on $X$ is integrable. 
\end{enumerate}
\end{theorem}

The quaternionic manifold satisfying \thetag{i} is called {\em hypercomplex}. 
Note that $\thetag{iii} \Rightarrow \thetag{ii} \Rightarrow \thetag{i}$ is 
obvious, so it suffices to prove $\thetag{i} \Rightarrow \thetag{iii}$. 

\section{Localization of $\protect\h$-modules.}\label{2}

\subsection{}
We begin with some linear algebra. 

\begin{defn}
A {\em quaternionic vector space} $V$ is a left module over the
algebra $\h$. 
\end{defn}

Let $V$ be a quaternionic vector space. Every algebra map $I:\C \to
\h$ defines by restriction a complex vector space structure on $V$,
which we will denote by $V_I$. Note that $\h$ is naturally a
left module over itself. The associated $2$-dimensional complex
vector space $\h_I$ is the same as in \ref{maps}, and we have $V_I =
\h_I \otimes_\h V$. 

\subsection{}\label{loc}
Let $\SB$ be the Severi-Brauer variety associated to the algebra
$\h$, that is, the variety of maximal right ideals in $\h$. By
definition $\SB$ is a real algebraic variety, an $\R$-twisted form
of the complex projective line $\cp$. It is also equipped with a
canonical maximal right ideal $\I \subset \h \otimes \calo$ in the
flat coherent algebra sheaf $\h \otimes \calo$ on $\SB$.

Let $V$ be a quaternionic vector space.  Consider the flat coherent
sheaf $V \otimes \calo$ on $\SB$ of right $\h \otimes \calo$-modules,
and let
$$
V_{loc} = V \otimes \calo / \I \cdot V \otimes \calo 
$$
be its quotient by the right ideal $\I$. Call the sheaf $V_{loc}$
{\em the localization} of the quaternionic vector space $V$.  The
localization is functorial in $V$ and gives a full embedding of the
category of quaternionic vector spaces into the category of flat
coherent sheaves on $\SB$.

Say that a flat coherent sheaf $\E$ on $\SB$ is {\em of weight $k$}
if the sheaf $\E \otimes \C$ on the complex projective line $\cp
\cong \SB \otimes \C$ is isomorphic to a sum of several copies of
the line bundle $\calo(k)$. The essential image of the localization
functor is the full subcategory of flat coherent sheaves of weight
$1$.

\subsection{}
The set $\cp \cong \SB(\C)$ of $\C$-valued points of the
variety $\SB$ was canonically identified in \ref{maps} with the set
$\Maps(\C,\h)$ of algebra maps $\C \to \h$. Let $I:\C \to \h$ be an
algebra map, and let $\wh{I} \in \cp$ be the corresponding
$\C$-valued point of the variety $\SB$ or, equivalently, the maximal
right ideal in the algebra $\h \otimes \C$. 

\begin{lemma}\label{hol}
Consider a quaternionic vector space $V$, and let $V_{loc}$ be its
localization. The fiber of the sheaf $V_{loc}$ at the point
$\wh{I} \in \cp$ is canonically isomorphic to the vector space
$V_I$. 
\end{lemma}

\proof Indeed, 
$$
V_I \cong \h_I \otimes_\h V \cong \h_I \otimes_{\h \otimes \C} V
\otimes \C \cong V \otimes \C / \wh{I} \cdot V \otimes \C =
V_{loc}|_{\wh{I}}, 
$$
and all the isomorphisms are canonical. 
\endproof 

\subsection{}
Consider now the set $\cp \cong \Maps(\C,\h)$ as the smooth
complex-analytic variety, and let $\V$ be the trivial bundle on
$\cp$ with the fiber $V \otimes_\R \C$. Since $\V$ is trivial, we
have a canonical holomorphic structure operator $\bar\6:\V \to
\A^{0,1}(\V)$ from $\V$ to the bundle $\A^{0,1}(\V)$ of $\V$-valued
$(0,1)$-forms on $\cp$.

The action of $\h$ on $V$ induces an operator $\I:\V \to \V$ which
acts as $I(\sqrt{-1})$ on the fiber $V$ of $\V$ at a point $I \in
\Maps(\C,\h)$. The operator $\I$ obviously depends smoothly on the
point $I$. It satisfies $I^2 = -1$ and induces therefore a smooth ``Hodge
type'' decomposition $\V = \V^{1,0} \oplus \V^{0,1}$. 

\begin{lemma}\label{locc}
The quotient $\V^{1,0}$ is compatible with the holomorphic structure
$\bar\6$ on $\V$. In other words, there exists a unique holomorphic
structure operator $\bar\6:\V^{1,0} \to \A^{0,1}(\V^{1,0})$ making
the diagram 
$$
\begin{CD}
\V            @>>>  \V^{1,0}         \\
@V{\bar\6}VV         @V{\bar\6}VV    \\
\A^{0,1}(\V)  @>>>  \A^{0,1}(\V^{1,0})
\end{CD}
$$
commutative. 
\end{lemma}

\proof This follows directly from Lemma~\ref{hol} by the usual
correspondence between flat coherent sheaves and holomorphic bundles
on the underlying complex-analytic variety. 
\endproof 

\section{Proof of the theorem.}\label{3}

\subsection{}
Let $Z$ be a smooth almost complex manifold. Let $\A^1(Z,\C)$ be the
complexified cotangent bundle to $Z$, and let $\A^\cdot(Z,\C)$ be
its exterior algebra. The almost complex structure on $Z$ induces
the Hodge type decomposition $\A^i(Z,\C) =
\oplus_{p+q=i}\A^{p,q}(Z)$.  Recall that the {\em Nijenhuis tensor}
$N$ of the almost complex manifold $Z$ is the composition
$$
N = P \circ d_Z \circ i:\A^{1,0}(Z) \to \A^1(Z,\C) \to
\A^2(Z,\C) \to \A^{0,2}(Z), 
$$
where $d_Z$ is the de Rham differential, $i:\A^{\cdot,0}(Z) \to
\A^\cdot(Z,\C)$ is the canonical embedding, and
$P:\A^\cdot(Z,\C) \to \A^{0,\cdot}(Z)$ is the canonical
projection. 

Recall also that the almost complex manifold $Z$ is called {\em
integrable} if its Nijenhuis tensor $N_Z:\A^{1,0}(Z) \to
\A^{0,2}(Z)$ vanishes. 

\subsection{}
We can now begin the proof of Theorem~\ref{main}. First we will
prove a sequence of preliminary lemmas. Let $M$ be a
smooth quaternionic manifold, and let $X$ be its twistor space. Since
by definition $X = M \times \cp$ as a smooth manifold, the cotangent
bundle $\A^1(X)$ decomposes canonically as 
\begin{equation}\label{eq.1}
\A^1(X) = \sigma^*\A^1(M) \oplus \pi^*\A^1(\cp), 
\end{equation}
where $\sigma:X \to M$, $\pi:X \to \cp$ are the canonical
projections. 

The almost complex structure $\I$ on $X$ preserves the decomposition
\eqref{eq.1}. Therefore \eqref{eq.1} induces decompositions 
\begin{align*}
\A^{1,0}(X) &= \A^{1,0}_M(X) \oplus \A^{1,0}_\cp(X), \\
\A^{0,2}(X) &= \A^{0,1}_M(X) \oplus \A^{0,1}_\cp(X),
\end{align*}
and, consequently, a decompositon 
\begin{equation}\label{dec}
\A^{0,2}(X) = \left( \A^{0,1}_M(X) \otimes \A^{0,1}_\cp(X) \right)
\oplus \A^{0,2}_M(X). 
\end{equation}
(Note that $\A^{0,1}_\cp(X)$ is of rank $1$, therefore
$\A^{0,2}_\cp(X)$ vanishes). More\-over, since the projection
$\pi:X \to \cp$ is compatible with the almost complex structures, we
have canonical isomorphisms 
$$
\A^{p,q}_\cp(X) \cong \pi^*\A^{p,q}(\cp). 
$$

\subsection{}
Let $N_X:\A^{1,0} \to \A^{0,2}(X)$ be the Nijenhuis tensor
of the almost complex manifold $X$. We begin with the following. 

\begin{lemma}
The restriction of the Nijenhuis tensor $N_X$ to the subbundle 
$$
\pi^*\A^{1,0}(\cp) \cong \A^{1,0}_\cp(X) \subset
\A^{1,0}(X)
$$ 
vanishes. 
\end{lemma}

\proof Indeed, since the map $\pi:X \to \cp$ is compatible with the
almost complex structures, the
diagram 
$$
\begin{CD}
\pi^*\A^{1,0}(\cp) @>>> \A^{1,0}(X)\\
@VVV                           @VV{N_X}V     \\
\pi^*\A^{0,2}(\cp) @>>> \A^{0,2}(X)
\end{CD}
$$
is commutative, and $\A^{0,2}(\cp)$ vanishes. 
\endproof 

Therefore the Nijenhuis tensor $N_X$ factors through a map
$$
N_X:\A^{1,0}_M(X) \to \A^{0,2}(X). 
$$

\subsection{}
Let now $N_X = N_1 + N_2$ be the decomposition of the Nijenhuis
tensor with respect to \eqref{dec}, so that $N_1$ is a map 
$$
N_1:\A^{1,0}_M(X) \to \A^{0,1}_M(X) \otimes \pi^*\A^{0,1}(\cp), 
$$
and $N_2$ is a map $N_2:\A^{1,0}_M(X) \to \A^{0,2}_M(X)$. 

\begin{lemma}
The component $N_1$ of the Nijenhuis tensor $N_X$ vanishes. 
\end{lemma}

\proof It suffices to prove that for every point $m \in M$ the
restriction $\wt{m}^*N_1$ of $N_1$ onto the corresponding twistor
line $\wt{m}:\cp \to X$ vanishes. Consider a point $m \in M$.  Let
$i:\wt{m}^*\A^{1,0}_M(X) \to \wt{m}^*\A^1_M(X,\C)$ be the canonical
embedding, and let
$$
P:\A^{0,1}(\cp) \otimes \wt{m}^*\A^1_M(X,\C) \to \A^{0,1}(\cp) \otimes
\wt{m}^*\A_M^{0,1}(X) 
$$
be the canonical projection. Since the twistor line $\wt{m}:\cp \to
X$ is compatible with the almost complex structures, we have
$\wt{m}^*\pi^*\A^{0,1}(\cp) \cong \A^{0,1}(\cp)$, and
\begin{multline*}
\wt{m}^*N_1 = P \circ \bar\6 \circ i:\wt{m}^*\A^{0,1}_M(X) \to
\wt{m}^*\A^1_M(X,\C) \to \\ 
\to \wt{m}^*\A^1_M(X,\C) \otimes \A^{0,1}(\cp) \to
\wt{m}^*\A^{0,1}_M(X) \otimes \A^{0,1}(\cp),
\end{multline*}
where $\bar\6:\wt{m}^*\A^1_M(X,\C) \to \wt{m}^*\A^1_M(X,\C) \otimes
\A^{0,1}(\cp)$ is the trivial holomorphic structure operator on the
constant bundle $\wt{m}^*\A^1_M(X,\C)$.

Let $V = T_mM$ be the tangent space to the manifold $M$ at the point
$m$. Since $M$ is quaternionic, $V$ is canonically a quaternionic
vector space. Let $\V$ and $\V^{1,0}$ be as in Lemma~\ref{locc}, and
let $\V^*$ and $(\V^{1,0})^*$ be the dual bundles on $\cp$. We have
canonical bundle isomorphisms
$$
\V^* \cong \wt{m}^*\A^1_M(X,\C) \qquad\qquad (\V^{1,0})^* \cong
\wt{m}^*\A^{1,0}_M(X) 
$$
compatible with the natural embeddings. By the statement dual to
Lemma~\ref{locc}, there exists a holomorphic structure operator
$\bar\6:\wt{m}^*\A^{1,0}_M(X) \to \wt{m}^*\A^{1,0}_M(X) \otimes
\A^{0,1}(\cp)$ making the diagram
$$
\begin{CD}
\wt{m}^*\A^{1,0}_M(X) @>{i}>> \wt{m}^*\A^1_M(X,\C) \\
@V{\bar\6}VV                    @V{\bar\6}VV     \\
\wt{m}^*\A^{1,0}_M(X) \otimes \A^{0,1}(\cp) @>{i \otimes \id}>>
\wt{m}^*\A^1_M(X,\C) \otimes \A^{0,1}(\cp) 
\end{CD}
$$
commutative. Therefore $N_1 = P \circ \bar\6 \circ i = P \circ (i
\otimes \id) \circ \bar\6$. But $P \circ (i \otimes \id) = 0$, hence
$N_1$ vanishes. 
\endproof 

\subsection{}
We can now prove Theorem~\ref{main}. As we have already proved, the
Nijenhuis tensor $N_X$ of the twistor space $X$ reduces to a bundle map 
$$
N_X:\A^{1,0}_M(X) \to \A^{0,2}_M(X). 
$$
This map vanishes identically if and only if for every point $m \in
M$ the restriction $\wt{m}^N_X$ of $N_X$ to the twistor line
$\wt{m}:\cp \to X$ vanishes. 

Consider a point $m \in M$. By Lemma~\ref{locc} the restriction
$\wt{m}^*\A^{1,0}_M(X)$ carries a natural holomorphic structure, and
it is a holomorphic bundle of weight $-1$ with respect to this
structure (in the sense of \ref{loc}). Consequently, the bundle
$\wt{m}^*\A^{0,2}_M(X)$ is a holomorphic bundle of weight
$2$. Moreover, the Nijenhuis tensor
$$
\wt{m}^*N_X = P \circ d \circ i:\wt{m}^*\A^{1,0}_M(X) \to
\A^{0,2}_M(X) 
$$
is a holomorphic bundle map. 

For every algebra map $I \in \Maps(\C,\h) \cong \cp$, the
restriction of the Nijenhuis tensor $N_X$ to a fiber $M \times I \in
M \times \cp = X$ of the projection $\pi:X \to \cp$ is the Nijenhuis
tensor for the induced almost complex structure $M_I$ on the
manifold $M$. Assume that Theorem~\ref{main}~\thetag{i} holds. Then
at least four distinct induced almost complex structures on $M$
corresponding to $I,J,\overline{I},\overline{J} \in \Maps(\C,\h)$
are integrable. Consequently, the Nijenhuis tensor $N_X$ vanishes
identically on fibers of the projection $\pi:X \to \cp$ over at
least four distinct points of $\cp$. Therefore the restriction
$\wt{m}^*N_X$ has at least four distinct zeroes. But as we have
proved, $\wt{m}^*N_X$ is a holomorphic map from a bundle of weight
$-1$ to a bundle of weight $2$. Therefore it vanishes
identically. Hence the almost complex manifold $X$ is integrable,
which finishes the proof of Theorem~\ref{main}.

\bigskip

\noindent
E-mail address: {\sc kaledin$@$balthi.dnttm.rssi.ru}

\end{document}